\documentclass{article}
\usepackage{amscd,amsfonts,euscript,amsmath,amsxtra,amssymb,amsthm}
\usepackage[english]{babel}
\newtheorem{theorem}{Theorem}

\newtheorem{proposition}{Proposition}

\theoremstyle{definition}

\newtheorem{remark}{Remark}
\theoremstyle{remark}

\numberwithin{equation}{section}

\usepackage[all]{xy}

\def\E{{{\mathbb E }}}

\def\EE{{{\mathcal E}}}

\def\OO{{{\mathcal O}}}

\def\Ext{{{\rm Ext }}}
\def\Tor{{{\rm Tor }}}

\def\hd{{{\rm hd \,}}}

\def\ker{{{\rm ker \,}}}

\def\red{{{\rm red}}}

\begin{document}

\medskip

\begin{center}
{\Large\sc A note on homological dimension of a family of coherent
sheaves}\end{center}
\medskip
\begin{center}
Nadezda V. TIMOFEEVA

\smallskip

Yaroslavl' State University

Sovetskaya str. 14, 150000 Yaroslavl', Russia

e-mail: {\it ntimofeeva@list.ru}
\end{center}
\bigskip

\begin{quote} We prove a theorem on how a conclusion on homological
dimension of the family of coherent sheaves  on a scheme can be
done from homological dimension of the restriction of this family
to the reduction of the base.
\\
Bibliography: 5 items.
\end{quote}

\markright{}
\begin{flushright}
{\it To the blessed memory of my Mum}
\end{flushright}
\medskip

The problem we solve in the present note is
\smallskip

{\it how to conclude about homological dimension of the family of
sheaves on the family of schemes if it is known about homological
dimension of the reduction?}
\smallskip

As usually if $T$ is a scheme with structure sheaf $\OO_T$ then
its {\it reduction} is a scheme consisting of the same topological
space but with structure sheaf equal to $\OO_{T_{\red}}:=
\OO_T/Nil(\OO_T)$ where $Nil (\OO_T)$ is nilradical of $\OO_T$.
This is called for brevity a {\it reduced scheme structure}. The
$\OO_T$-module epimorphism onto the quotient module sheaf
$\OO_T\twoheadrightarrow \OO_{T_{\red}}$ induces a canonical
closed immersion $T_{\red}\hookrightarrow T$ of schemes. From now
$T$ is a base scheme of a flat morphism of finite type $f: X\to T$
of Noetherian schemes. We introduce  notations $X_{\red}:= X
\times _T T_{\red}$, $f_{\red}: X_{\red} \to T_{\red}$, i.e.
$X_{\red}$ is a restriction of the family $X$ to the reduction
$T_{\red}$ as to a closed subscheme in $T$. Actually the scheme
$X_{\red}$ can be nonreduced but $f_{\red}$ is flat morphism as a
morphism obtained by a base change of a flat morphism. Now let
$\E$ be a coherent $\OO_X$-module and let $\E$ be flat as
$\OO_T$-module. This is a standard situation in various problems
of algebraic geometry when families of sheaves are considered,
especially in moduli problems. Denote $\E_{\red}=\E
\otimes_{\OO_T} \OO_{T_{\red}}$.

Let $A$ be a commutative ring, $M$ $A$-module. Also introduce
parallel algebraic notations: a {\it reduction} $A_{\red}:=A/Nil
A$ of the ring $A$ is its quotient ring over its nilradical. If a
commutative ring $B$ is an $A$-algebra then we denote $B\otimes _A
A_{\red}=: B_{\red}$. Hence $B_{\red}$ {\it is not obliged to be
reduced} but it is $A_{\red}$-flat whenever $B$ is $A$-flat, by
well-known change-of-ring theorem. Also denote $M\otimes _A
A_{\red}=:M_{\red}.$

We say that $M$ has {\it homological dimension not greater then}
$n$ (notation: $\hd_A M \le n$) if one of following equivalent
conditions holds
\cite[Ch. 7, theorem 1.1]{MacLane}: \\
(1) for all $A$-modules $N$ $\Ext^{n+1}_A(M,N)=0$, \\
(2) in any exact $A$-sequence $0\to E_n \to E_{n-1} \to \dots \to
E_0 \to M \to 0$ if $E_j$ are projective for $0\le j\le n-1$ then
$E_n$ is also projective,\\
(3) there is a projective $A$-resolution of length $n$, i.e. there
exists an exact sequence
$$
0\to E_n \to \dots \to E_1 \to E_0 \to M \to 0
$$
with $E_j$ projective for $0\le j  \le n$.

Since if $A$ is local ring then any projective $A$-module is free
\cite[theorem 19.2 and comment thereafter]{Eisen}. Then when
working with coherent sheaves on schemes we speak of locally free
resolutions instead of projective ones.

As usually the symbol $\hd_X \E$ means homological dimension  of
the coherent sheaf $\E$ as $\OO_X$-module.

We prove the following well-expected result. \begin{theorem} {\bf
Algebraic version.} Let $f^{\sharp}:A\to B$ be local homomorphism
of local Noetherian rings and $B$ is flat as $A$-algebra. Let $M$
is $B$-module of finite type which is flat over $A$.  Then
following assertions are equivalent:\\ 1) $\hd_B M \le n$, \\
2) $\hd_{B_{\red}} M_{\red} \le n$.
\end{theorem}
\begin{theorem} {\bf Sheaf version.} Let $f: X\to T$ be a flat
morphism of Noetherian schemes, $\E$ coherent $\OO_X$-module which
is flat over $T$. Then following assertions are equivalent:\\
1) $\hd_{X}\E \le n$, \\
2) $\hd_{X_{\red}}\E_{\red} \le n$.
\end{theorem}
\begin{remark} The case $n=1$ for a trivial family of schemes over a
field was considered in the author's paper
\cite{Tim8}.
\end{remark}
\begin{remark} Since both theorems are just versions of the same
result their proofs are transferred literally to each other and we
prove an algebraic version.
\end{remark}
\begin{proof}[Proof of Theorem 1]  For the implication
1)$\Rightarrow$2) we do not need locality and  consider $B$-exact
sequence
$$0\to E_n \to E_{n-1} \to \dots \to E_1 \to E_0 \to M \to 0
$$
with $E_j$ free for $j\ge 0$. Cutting it into triples we have
\begin{eqnarray*}
0\to E_n \to E_{n-1} \to M^{(n-1)} \to 0,\\
\dot \quad\dot \quad\dot \quad\dot \quad\dot \quad\dot \quad\dot
\quad\dot \quad
\dot \quad\dot \quad\dot \quad\dot \quad\dot \quad\dot \quad,\\
0\to M^{(j+1)} \to E_j \to M^{(j)} \to 0,\\
\dot \quad\dot \quad\dot \quad\dot \quad\dot \quad\dot \quad\dot
\quad\dot \quad
\dot \quad\dot \quad\dot \quad\dot \quad\dot \quad\dot \quad,\\
0\to M^{(1)} \to E_0 \to M \to 0.
\end{eqnarray*}
Since $M$ and $E_j$, $j\ge 0$, are flat $A$-modules \cite[Ch. 3,
sect. 7, transitivity (1)]{Mats} we have that $M^{(j)}$ are
$A$-flat for $j\ge 1$. Hence tensoring by $\otimes_B
B_{\red}=\otimes _A A_{\red}$ ("cancelation formula") we have
$\Tor_i^A(M^{(j)}, A_{\red})=\Tor_i^A(M, A_{\red})=0$ and come to
exact sequence
$$0\to E_{n\red} \to E_{(n-1)\red} \to \dots \to E_{1\red} \to
E_{0\red} \to M_{\red} \to 0
$$ with $E_{j\red}$ free as $B_{\red}$-modules for $0\le j \le n$.

For the opposite implication we organize induction on $n$. Assume
that the theorem is true for coherent $B$-sheaves of homological
dimension not greater then $n-1$. Let we are given an $B$-module
epimorphism $E_0\twoheadrightarrow M$ where $E_0$ is free and let
$M':=\ker (E_0\twoheadrightarrow M)$. Let
$\hd_{B_{\red}}M_{\red}\le n$. We are to conclude that $\hd_B M
\le n$. By the exact $B$-triple
$$0\to M' \to E_0 \to M \to 0
$$
and tensoring by $\otimes_B B_{\red}=\otimes_A A_{\red}$,
since $M$ is $A$-flat then  $\Tor^A_i(M, A_{\red})=0$ for $i>0$
and hence
 we come to the exact $A_{\red}$-triple
$$0\to M'_{\red} \to E_{0\red} \to M_{\red} \to 0
$$
where $\hd_{A_{\red}}M'_{\red} \le n-1$. By flatness of $E_0$ as
$B$-module and of the homomorphism $f^{\sharp}$ the term $E_0$ is
$A$-flat \cite[Ch. 3, sect. 7, transitivity (1)]{Mats} and hence
$M'$ is also flat as $A$-module. By the inductive assumption
$\hd_B M'\le n-1$ and hence $\hd_B M \le n$.

For the base of induction set $n=0$. This means that $M$ is flat
as $A$-module  and $M_{\red}$ is  free as $A_{\red}$-module.

Apply the  following result from A. Grothendieck's SGA:
\begin{proposition}\cite[Ch. IV, Corollaire 5.9]{SGAI} Let $A\to B \to
C$ be local homomorphisms of local Noetherian rings, $M$ be
$C$-module of finite type. Assume that $B$ is flat over $A$ and
$k$ is a residue field of $A$. Then following assertions are
equivalent:\\
(i) $M$ is $B$-flat;\\
(ii) $M$ is $A$-flat and $M\otimes _A k$ is $B\otimes_A k$-flat.
\end{proposition}

For our purposes set $B\to C$ to be an identity isomorphism, $M$
is flat over $A$ and $M_{\red}$ is free (i.e. flat) as
$A_{\red}$-module. Then   $M \otimes_A k=M_{\red} \otimes
_{A_{\red}} k$, $B\otimes _A k=B_{\red}\otimes_{A_{\red}} k$, and
$M_{\red} \otimes_{A_{\red}} k$ is flat over
$B_{\red}\otimes_{A_{\red}} k$ because $M_{\red}$ is free over
$B_{\red}$. From this we conclude that $M$ is  free as $B$-module.
This completes the proof of the theorem 1.
\end{proof}

\end{document}